\documentclass[12pt]{article}
\begin{document}
\bibliographystyle{alpha}
\bibliography{}

\begin{thebibliography}{BGS}

\bibitem[Ba]{Ba} Baragar, A., ``Rational points on $K3$ surfaces in
$\P^1\times\P^1\times\P^1$'', Math. Ann. 305, no. 3, 541--558 (1996)

\bibitem[BGS]{BGS} Bost, J.-B., Gillet, H., and Soul\'e, C., ``Heights
of Projective Varieties and Positive Green Forms'', J. Amer.
Math. Soc. 7, no. 4, 903--1027 (1994)

\bibitem[Bi]{Bi} Billard, H., ``Propri\'et\'es arithm\'etiques d'une
famille de surfaces $K3$'', Compositio Math. 108, no. 3, 247--275
(1997)

\bibitem[BM]{BM} Batyrev, V. and Manin, Yu., ``Sur le nombre de points
rationnels de hauteur born\'ee des vari\'et\'es alg\'ebriques'',
Math. Ann. 286, 27-43 (1990)

\bibitem[IS]{IS} Iskovskikh, V.A. and Shafarevich, I.R., ``Algebraic
Surfaces'', in {\bf Algebraic Geometry II}, Encyclopedia of
Mathematical Sciences, vol. 35, Springer-Verlag, Berlin, 1996.

\bibitem[KT]{KT} King, H. and Todorov, A., ``Rational points on some
Kummer surfaces'', preprint (1992)

\bibitem[Sch]{Sch} Schanuel, ``Heights in number fields'',
Bull. Soc. Math. France 107 433-449 (1979)

\bibitem[SD]{SD} Saint-Donat, B.  ``Projective Models of $K3$
Surfaces'', Amer. J. Math., Vol. 96, No. 4, 602--639 (1974)

\bibitem[Si]{Si} Silverman, J., ``Computing heights on $K3$ surfaces:
a new canonical height'', Invent. Math. 105, 347-373 (1991)

\bibitem[Ts]{Ts} Tschinkel, Yu., Ph.D. thesis, MIT (1992)

\end{thebibliography}
\newtheorem{prop}{Proposition}
\newtheorem{cor}{Corollary}
\newtheorem{lem}{Lemma}
\newtheorem{thm}{Theorem}
\newtheorem{dfn}{Definition}
\def\O{{\mathcal O}}
\def\F{{\mathcal F}}
\def\N{{\mathbf N}}
\def\Q{{\mathbf Q}}
\def\R{{\mathbf R}}
\def\C{{\mathbf C}}
\def\P{{\mathbf P}}
\def\Z{{\mathbf Z}}
\def\E{{\exists}}
\def\A{{\forall}}
\def\M{{\mathcal M}}
\def\tr{\mbox{Tr}}
\def\qed{\tiny $\clubsuit$ \normalsize}
\def\qedl{\tiny $\clubsuit$ \normalsize}
\title{Counting Rational Points on K3 Surfaces}
\author{David McKinnon}
\maketitle

\begin{abstract}
For any algebraic variety $V$ defined over a number field $k$, and
ample height function $H_D$ on $V$, one can define the counting function
$N_{V,D}(B) = \#\{P\in V(k) \mid H_D(P)\leq B\}$.  In this paper, we calculate
the counting function for Kummer surfaces $V$ whose associated abelian
surface is the product of elliptic curves.  In particular, we
effectively construct a finite union $C = \cup C_i$ of curves
$C_i\subset V$ such that $N_{V-C,D}(B)\ll N_{C,D}(B)$; that is, $C$ is an
accumulating subset of $V$.  In the terminology of Batyrev and Manin
[BM], this amounts to proving that $C$ is the first layer of the
arithmetic stratification of $V$.
\end{abstract}

\section{Introduction}

Counting rational points on algebraic varieties is one of the
fundamental questions of number theory.  However, if an algebraic
variety contains infinitely many rational points one must define the
question more precisely.  The most natural way to do this is to define
a notion of density on the set of rational points.  This density is
calculated with respect to a height $H$, which assigns a real number
to a rational point $P$.  Thus, for a variety $V$ defined over a
number field $K$ and an ample divisor $D$ on $V$, we study the
counting function:
\[N_{V,D}(B) = \mathrm{card}\{P\in V(K) \mid H_D(P)\leq B\}\]
and investigate the properties of $N_{V,D}(B)$ as $B$ gets arbitrarily
large.  This function may be radically different for different choices
of $D$.  Also, since $H_D$ is only defined up to multiplication by a
bounded function, so too is the definition of $N_{V,D}$.

In the case of $K3$ surfaces, this question has been investigated by
many people.  Silverman [Si] introduced a canonical height on $K3$
surfaces embedded in $\P^2\times\P^2$, analogous to the canonical
height on an elliptic curve.  Baragar [Ba] extended Silverman's
results to other $K3$ surfaces.  Although both authors obtain theorems
about the distribution of rational points in orbits of certain group
actions, neither was able to obtain estimates of the global counting
function.  Billard [Bi] has recently extended their results still
further, and gives an estimate for $N_{V,D}(B)$ in a certain case.

Another approach was taken by Tschinkel [Ts], who develops a theory of
finite heights to obtain estimates of $N_{V,D}(B)$ for some rational
surfaces, and upper bounds on $N_{V,D}(B)$ for some $K3$ and Enriques
surfaces.  King and Todorov [KT] use the results of [Ts] to estimate
$N_{V,D}(B)$ for a certain class of Kummer surfaces admitting a double
cover of a del Pezzo surface.

In this paper, the particular $K3$ surfaces we will study are Kummer
surfaces $V$ whose associated abelian surfaces are isomorphic to a
product of elliptic curves.  We define an 18-dimensional cone
${\mathcal C}$ of ample divisors in the N\'eron-Severi lattice of $V$,
and calculate the value of $N_{V,D}(B)$ with respect to an arbitrary
divisor in ${\mathcal C}$.  More specifically, if we measure heights
with respect to a divisor $D\in{\mathcal C}$, we will show that
$N_{V,D}(B)$ is asymptotically equal to $N_{C,D}(B)$, where $C$ is the union
of all rational curves of minimal $D$-degree on $V$.  We also
calculate explicitly which curves lie in $C$.

Batyrev and Manin [BM] have introduced a refinement of the counting
function called the {\it arithmetic stratification}.  Roughly
speaking, a subset $W$ of $V$ is said to be {\it accumulating} with
respect to an ample divisor $D$ if most of the rational points of $V$
lie on $W$, where heights are measured with respect to $D$.  That is,
if $\lim_{B\rightarrow \infty} N_{V,D}(B)/N_{W,D}(B) = 1$.  The arithmetic
stratification of a variety $V$ is an ascending chain of Zariski
closed subsets $W_1\subset W_2\subset W_3\subset\ldots$ with the
property that $W_i - W_{i-1}$ is an accumulating subset of $V -
W_{i-1}$.  $W_i$ is said to be the $i$th layer of the arithmetic
stratification.  Since layers in the arithmetic stratification are
typically finite unions of rational curves, the value of $N_{V,D}(B)$ will
immediately follow from Schanuel's theorem [Sch], which calculates the
counting function for $\P^n$.

Corollary 1 explicitly identifies the first layer of the arithmetic
stratification of $V$ with respect to an 18-dimensional cone of ample
divisors.  The number of rational points lying on any given rational
curve on $V$ can be easily calculated from Schanuel's theorem; the
hard part comes from Theorem 1, which estimates the number of rational
points on the complement of the union of these curves.  By comparing
the counting functions for certain rational curves constructed on $V$
with the counting function for the complement $U$ of the union of
these curves, the structure of the top layers of the arithmetic
stratification is revealed.

I would like to thank Jim Bryan, Yuri Tschinkel, and Tom Tucker for
helpful comments and conversations.  I am especially grateful to Paul
Vojta, without whose advice and support this paper would never have 
appeared.

\section{Geometric Preliminaries}

Let $C_1$ and $C_2$ be elliptic curves defined over some number field
$K$, such that all points of order 1 and 2 on the curve are also
defined over $K$.  Let $A$ be the product $C_1\times C_2$.  Let
$i:A\rightarrow A$ be the involution $i(x,y) = (-x,-y)$, and let
$S$ be the quotient of $A$ by $i$.  Then there is a 2-to-1 map
$q:A\rightarrow S$ which is ramified at 16 points; namely, the
points $(a,b)$, where $a$ and $b$ are points of order 1 or 2.  It
turns out that these 16 points are rational double points of
$S$, which is smooth away from them.  

By blowing up these 16 points, one constructs a smooth surface
$p:X\rightarrow S$, which is a $K3$ surface defined over $K$ [IS].
This construction can be done with an arbitrary abelian surface $A$,
and the resulting $K3$ surface is called the Kummer surface associated
to the abelian surface $A$.  

Let $\pi_i:A\rightarrow C_i$ be the projection maps, and let
$F^{\prime}_i$ be the algebraic equivalence class of fibres of
$\pi_i$.  This induces a pair of algebraic equivalence classes $F_i =
p^*q_* F^{\prime}_i$ on $X$ -- since algebraic and linear equivalence
are identical on a $K3$ surface [IS], these are divisor classes on
$X$.  Thus, the maps $\pi_i$ descend to maps $\tilde{\pi}_i:
X\rightarrow \P^1$.

Denote the 16 singular points of $S$ by $(a_i,b_j)$, $1\leq i,j\leq
4$, where $a_i$ and $b_j$ denote the 2-division points on $C_1$ and
$C_2$, respectively, and let $E_{ij}$ denote the corresponding
exceptional divisors on $S$.  For each $i$, $1\leq i \leq 4$, the
divisor $B_i = p^*q_*\pi_1^*{a_i}$ is the union of the four curves
$E_{ij}$, $1\leq j \leq 4$, and the strict transform of
$q_*(\{a_i\}\times E_2)$.  By the theory of singular fibres of
elliptic surfaces [IS], it follows that this strict transform is a
double curve, which is smooth and rational in its induced reduced
structure.  Thus, we may write $F_1 \equiv B_i = \sum_{j=1}^4 E_{ij} +
2L_i$, where $L_i$ is a smooth rational curve.  Similarly, we may
write $F_2 \equiv \sum_{i=1}^4 E_{ij} + 2M_j$, where $M_j$ is a smooth
rational curve.

Using the adjunction formula and elementary properties of intersection
theory, it is not hard to verify the following intersection numbers:
\[\begin{array}{ll}
L_i^2 = M_i^2 = E_{ij}^2 = -2 & L_iM_j = 0 \\
F_1L_i = F_2M_i = 0 & F_1M_i = F_2L_i = 1 \\
L_iL_j = M_iM_j = 0 \hspace{.5in} \mbox{(if $i\neq j$)} &
\end{array}\]
Let $S$ and $T$ be non-empty subsets of $N_4=\{1,2,3,4\}$.  Define
divisors \[A_{S,T} = (\mathrm{card}(S))F_1 + (\mathrm{card}(T))F_2 -
\sum_{i\in S,j\in T} E_{ij}\] These divisors, together with $F_1$ and
$F_2$, span a rank 18 sublattice of Pic($X$), and therefore an
18-dimensional subspace of the vector space $NS_\R(X) =
\mathrm{Pic}(X)\otimes\R$.  For a generic choice of $C_1$ and $C_2$,
$NS_\R(X)$ has dimension 18 [IS], so the divisors $A_{S,T}$, and $F_i$
span all of $NS_\R(X)$ for such $X$.

Moreover, for any ample divisor $D$, write $D=d_1F_1 + d_2F_2 + \sum
e_{ij}E_{ij}$.  Since $D.E_{ij}>0$ and $E_{ij}^2=-2$, we must have 
$e_{ij}<0$.  Therefore, it follows that any ample divisor $D$ on $X$ can
be written as 
\begin{equation}
D = \sum_{S,T} a_{S,T}A_{S,T} + c_1F_1 + c_2F_2
\end{equation}
\noindent
where $a_{S,T}\geq 0$.  Also note that if $D$ is written in this form,
then we may assume without loss of generality that $a_{N_4,N_4} =
\min\{e_{ij}\}$.  This will be assumed to be true in all that
follows.

\section{The Main Theorem}

We are now ready to state the main theorem.  All counting functions are
defined with respect to the height associated to the divisor $D$.

\begin{thm}
Let $D$ be an ample divisor on $X$ written as in (1).  Assume that
$a_{S,T}$ are non-negative rational numbers, and $c_i$ are rational
numbers.  Define:
\[\gamma_1 = \sum_{S,T} \mathrm{card}(S)a_{S,T}, \hspace{.75in} \gamma_2 = 
\sum_{S,T} \mathrm{card}(T)a_{S,T}\]
\[\alpha = \max\{\frac{2\gamma_1 + 2\gamma_2}{\gamma_1\gamma_2 + 
\gamma_2c_1 + \gamma_1c_2}, \frac{2}{\gamma_1 + c_1}, 
\frac{2}{\gamma_2 + c_2}\}\]
\noindent
Define $U = V\setminus \bigcup R$, where $R$ ranges over all smooth
rational curves on $V$ of the form $E_{ij}$, $L_i$, or $M_i$.  Assume
that $\gamma_1\gamma_2 + \gamma_2c_1 + \gamma_1c_2 > 0$.  Then:

\vspace{.1in}

i) If $\alpha = \frac{2\gamma_1 + 2\gamma_2}{\gamma_1\gamma_2 +
\gamma_2c_1 + \gamma_1c_2}$ and either $c_1 = \gamma_2 + c_2$ or $c_2
= \gamma_1 + c_1$, then $N_{U,D}(B) = O(B^{\alpha}\log B)$.

\vspace{.1in}

ii) If $\alpha = \frac{2}{\gamma_1 + c_1}$ and $c_2 = \gamma_1 + c_1$,
then $N_{U,D}(B) = O(B^{\alpha}\log B)$.

\vspace{.1in}

iii) If $\alpha = \frac{2}{\gamma_2 + c_2}$ and $c_1 = \gamma_2 + c_2$,
then $N_{U,D}(B) = O(B^{\alpha}\log B)$.

\vspace{.1in}

iv) If none of the previous three cases occur, then $N_{U,D}(B) =
O(B^\alpha)$.

\end{thm}

\begin{cor}
Let $X$ be a $K3$ surface as described above, and let $D$ be an ample
divisor, written as in Theorem 1.  Write $A =
\min\{D.E_{ij},D.L_i,D.M_i\}$, and assume that the following
inequality holds:
\begin{equation}
A(\gamma_1+\gamma_2) < \gamma_1\gamma_2 + \gamma_2c_1 + \gamma_1c_2
\label{cond}
\end{equation}
Then the counting function for $X(K)$ is given by $N_{X,D}(B) = cB^{8/A} +
E(B)$, where $c$ is a constant depending only on $K$, $X$, and the
choice of height function $H_D$, and $E(B)=O(B^q)$ is an error term
with an easily calculable $q<8/A$.  Moreover, the main term measures
the number of rational points lying on the union of rational curves of
minimal $D$-degree, which must all be of the form $E_{ij}$, $L_i$, or
$M_i$, and $E(B)$ bounds the number of rational points not lying on
such curves.
\end{cor}

More precisely, we prove that the first term in the above expression
represents the number of rational points lying on smooth rational
curves of minimal $D$-degree on $X$.  The error term represents the
combination of Schanuel's error term, the estimate from Theorem 1 for
$N_{U,D}(B)$, and the number of rational points lying on the curves
$E_{ij}, L_i,$ and $M_i$ of non-minimal $D$-degree.  

Put another way, the first layer of the arithmetic stratification (as
defined by Batyrev and Manin [BM]) of $X$ with respect to an ample
divisor $D$ is the union of all smooth rational curves of minimal
$D$-degree, provided that $D$ can be expressed in a form for which
inequality~$(\ref{cond})$ is satisfied.  Moreover, these curves are
all of the form $E_{ij}$, $L_i$, or $M_i$.

\vspace{.1in}

\noindent
{\it Proof of Corollary 1:} \/ The following curves have the following
degrees:
\[\deg_D(E_{mn}) = 2\sum_{S\ni m,T\ni n}a_{S,T} \]
\[\deg_D(L_n) = c_2 + \sum_{S\not\ni n}\sum_{T} \mathrm{card}(T) a_{S,T}\]
\[\deg_D(M_n) = c_1 + \sum_{S}\sum_{T\not\ni n} \mathrm{card}(S) a_{S,T}\] 
\noindent
and by Schanuel's Theorem for a smooth rational curve $C$ of degree
$d$ in projective space, we have $N_{C,\O(d)}(B) = cB^{2/d} + O(B^{2/d -
1/Nd})$, where $N = [K:\Q]>1$ and $c$ is a complicated constant,
calculated explicitly by Schanuel.  (In the special case $K=\Q$, the
error term must be replaced by $O(B^{1/d}\log B)$.)

It suffices to show that $N_{U,D}(B) < N_{E_{ij},D}(B)$ for sufficiently
high $B$, where $E_{ij}$ is the exceptional curve of lowest degree.
By Schanuel's Theorem, we have $N_{E_{ij},D}(B) = cB^{2/A} + O(B^{2/A -
2/NA})$, where $c$ is a constant depending only on $K$, $X$, and $D$,
and $N=[K:\Q]$.  (If $K=\Q$, the error term must be appropriately
modified.)  By the theorem, then, it suffices to show that
$\frac{2}{A} > \alpha$.

If $\alpha = \frac{2\gamma_1 + 2\gamma_2}{\gamma_1\gamma_2 + 
\gamma_2c_1 + \gamma_1c_2}$, then the desired inequality follows
immediately from equation (\ref{cond}).

Assume that $\alpha = \frac{2}{c_1 + \gamma_1}$.  Since $\deg_D M_n >
0$ for $n=1,2,3,4$, it follows that:
\begin{eqnarray*}
0 & < & c_1 + \sum_{S}\sum_{T\not\ni n} \mathrm{card}(S) a_{S,T} \\ 
& < & c_1 + A + \sum_{S\neq N_4\neq T} \mathrm{card}(S) a_{S,T} \\
& < & c_1 - A + \gamma_1
\end{eqnarray*}
\noindent
which implies immediately that $\frac{2}{A} > \alpha$, as desired.
Similarly, if $\alpha = \frac{2}{c_2 + \gamma_2}$, then $\frac{2}{A} >
\alpha$ follows from the positivity of $\deg L_n$.

Finally, we must establish that there are no rational curves of minimal
$D$-degree other than those of the form $E_{ij}$, $L_{i}$, or $M_i$.
Assume there exists such a curve $C$ of minimal $D$-degree.  Then
$C\cap U$ is a dense open subset of $U$, so $N_{C,D}(B)\ll N_{U,D}(B)$.  But
for a curve of minimal $D$-degree, we have just established that
$N_{C,D}(B)\not\ll N_{U,D}(B)$.  The corollary follows.  \qed

\vspace{.1in}

\noindent
{\it Proof of Theorem 1:} \/ The key idea is to estimate an
arbitrary height function $H_L$ in terms of the height functions
$H_{F_i}$, which are easily computed.  The first step is to note that
the divisors $A_{S,T}$, $F_1$, and $F_2$ span a rank 18 sublattice of
Pic($X$).  This can be proven by explicit calculation.  For
non-isogenous elliptic curves $C_1$ and $C_2$ (as is generally the
case), Pic($X$) is a free $\Z$-module of rank 18.  Therefore, height
calculations with respect to a general ample sheaf $L$ can be reduced
to calculations with respect to the divisors $A_{S,T}$, $F_1$, and
$F_2$.  Write $H_1 = H_{F_1}$ and $H_2 = H_{F_2}$.

\begin{lem}
Let $U=X\setminus (\cup E_{ij} \cup M_i \cup L_i)$.  Let $D=A_{S,T}$.
Then for any point $P\in U(K)$:
\[\max\{H_1(P)^{\mathrm{card}(S)}, H_2(P)^{\mathrm{card}(T)}\} 
\ll H_D(P)\]
\[H_D(P) \ll H_1(P)^{\mathrm{card}(S)}H_2(P)^{\mathrm{card}(T)}\]
\end{lem}

\noindent
{\it Proof:} \/ The second inequality follows immediately from the
arithmetic Bezout theorem of [BGS] (or from the effectivity of the
divisors $E_{ij}$), so it suffices to prove the first inequality.
This follows immediately from the effectivity of the divisors $D -
\mathrm{card}(S)F_1$ and $D - \mathrm{card}(T)F_2$.  \qed

\vspace{.1in}

From Lemma 1, we get the following inequalities:
\[H_1(P)^{c_1}H_2(P)^{c_2} \prod_{S,T}\max\{H_1(P)^{\mathrm{card}(S)},
H_2(P)^{\mathrm{card}(T)}\}^{a_{S,T}} \ll H_D(P)\]
\[H_D(P) \ll H_1(P)^{\gamma_1+c_1}H_2(P)^{\gamma_2+c_2}\]
These estimates are enough to prove Theorem 1.  The first inequality
above implies for any point $P\in U(K)$ (since $a_{S,T}\geq 0$):
\[H_D(P) \gg \max\{H_1(P)^{\gamma_1+c_1}H_2(P)^{c_2}, 
H_1(P)^{c_1}H_2(P)^{\gamma_2+c_2}\}\]

The number of points of points of height at most $B$ on $U$ is
therefore bounded by the number of integer lattice points contained in
a certain plane region times a constant factor (which is immaterial to
the result of the theorem).  By Schanuel's Theorem, there are $\ll B^2$
points of height at most $B$ on $\P^1$ with respect to the height
attached to $\O(1)$, so that if $H_1(P) \leq B$, then there are
$\ll B^2$ choices for $\tilde\pi_1(P)$, and similarly for $H_2(P)$.
Therefore, set $x=H_1(P)^2$ and $y=H_2(P)^2$.  If $H_D(P) \leq B$,
then we get
\[\max\{x^{\gamma_1+c_1}y^{c_2},x^{c_1}y^{\gamma_2+c_2}\} \leq B^2\]
Thus, $N_{U,D}(B)$ is bounded above by a constant factor times the
number of lattice points contained in the plane region $R$ defined by the
inequalities:
\[R = \{(x,y)\in\R^2 \mid \mbox{$x\geq 1$, $y\geq 1$, 
$x^{\gamma_1+c_1}y^{c_2}\leq B^2$, {\rm and} 
$x^{c_1}y^{\gamma_2+c_2}\leq B^2$}\}\]  
This is asymptotically equal to the area of this region (again, up to
an irrelevant constant factor), plus two extra terms counting lattice
points lying on the boundary lines $x=1$ and $y=1$.  This may be
computed as follows.

\vspace{.1in}

\noindent
CASE I: $c_2>0$.  Define $\delta = 2(\gamma_1\gamma_2 +
\gamma_2c_1 + \gamma_1c_2)^{-1}$.  The two curves
$x^{\gamma_1+c_1}y^{c_2}=B^2$ and $x^{c_1}y^{\gamma_2+c_2}=B^2$
intersect at the point $(B^{\delta\gamma_2},B^{\delta\gamma_1})$.  Thus,
the number of lattice points inside $R$ may be computed by:
\begin{eqnarray*}
& & \int_1^{B^{\delta\gamma_2}}
(B^2x^{-c_1})^{\frac{1}{\gamma_2+c_2}} dx + \int_{B^{\delta\gamma_2}}
^{B^{2/\gamma_1+c_1}}(B^2x^{-\gamma_1-c_1})^{\frac{1}{c_2}} dx 
+ B^{\frac{2}{\gamma_1+c_1}} + B^{\frac{2}{\gamma_2 + c_2}} \\
& = & -(\frac{\gamma_2 + c_2}{\gamma_2 + c_2 -
c_1})(B^{\frac{2}{\gamma_2+c_2}} - B^{\delta(\gamma_1 + \gamma_2)}) \\
& & -(\frac{c_2}{\gamma_1+c_1-c_2})(B^{\frac{2}{\gamma_1+c_1}} -
B^{\delta(\gamma_1 + \gamma_2)}) + B^{\frac{2}{\gamma_1+c_1}} + 
B^{\frac{2}{\gamma_2 + c_2}} \\
& = & O(B^\alpha)
\end{eqnarray*}
\noindent
unless $\gamma_1+c_1=c_2$ or $\gamma_2+c_2=c_1$, in which case obvious
modifications to the computation will give the desired result.  (Note
that $\gamma_i+c_i>0$ for all $i$ by the positivity of $\deg_D L_n$
and $\deg_D M_n$.)

\vspace{.1in}

\noindent
CASE II: $c_2 = 0$.  Retain the notation of the previous case.  The number
of lattice points lying inside $R$ is now bounded by:
\begin{eqnarray*}
& & \int_1^{B^{\frac{2}{\gamma_1+c_1}}}
(B^2x^{-c_1})^{\frac{1}{\gamma_2}} dx + B^{\frac{2}{\gamma_1+c_1}} +
B^{\frac{2}{\gamma_2}} \\ 
& = & -(\frac{\gamma_2}{\gamma_2-c_1})(B^{\frac{2}{\gamma_2}}
- B^{\frac{2}{\gamma_2+c_1}}) + B^{\frac{2}{\gamma_1+c_1}}
+ B^{\frac{2}{\gamma_2}} \\
& = & O(B^\alpha)
\end{eqnarray*}
\noindent
again with the obvious modifications in the case that $c_1=\gamma_2$.

\vspace{.1in}

\noindent
CASE III: $c_2 < 0$.  Again retaining the notation of the previous cases,
we may compute the number of lattice points lying inside $R$ by:
\begin{eqnarray*}
& & \int_1^{B^{\delta\gamma_2}}
(B^2x^{-c_1})^{\frac{1}{\gamma_2+c_2}} 
- \max\{(B^2x^{-\gamma_1-c_1})^{\frac{1}{c_2}},1\} dx 
+ B^{\gamma_2\delta} + B^{\frac{2}{\gamma_2+c_2}} \\
& \leq & \int_1^{B^{\delta\gamma_2}}
(B^2x^{-c_1})^{\frac{1}{\gamma_2+c_2}} dx
+ B^{\delta\gamma_2} + B^{\frac{2}{\gamma_2+c_2}} \\
& = & -(\frac{\gamma_2 + c_2}{\gamma_2 + c_2 -
c_1})(B^{\frac{2}{\gamma_2+c_2}} - B^{\delta(\gamma_1 + \gamma_2)})
+ B^{\frac{2}{\gamma_1+c_1}} + B^{\frac{2}{\gamma_2 + c_2}} \\
& = & O(B^\alpha)
\end{eqnarray*}
\noindent
again with obvious modifications in the case that $c_1=c_2+\gamma_2$.

Thus, the proof of Theorem 1 is complete.  \qed

\noindent
David McKinnon \\
Department of Mathematics \\
University of California at Berkeley \\
Berkeley, CA 94720 \\
E-mail: mckinnon@math.berkeley.edu

\end{document}